\documentclass[a4paper,12pt]{article}

\input{mymacros.sty}

\newcommand{\SR}{{SR}}

\newcommand{\bSigma}{{\boldsymbol{\Sigma}}}
\newcommand{\bSigmabar}{{\overline{\bSigma}}}
\newcommand{\bfDbar}{{\overline{\bfD}}}
\newcommand{\scHom}{\mathop{\scH o m}\nolimits}

\newcommand{\bvarphi}{{\boldsymbol{\varphi}}}
\newcommand{\bscE}{{\boldsymbol{\scE}}}
\newcommand{\bscL}{{\boldsymbol{\scL}}}
\newcommand{\bscV}{{\boldsymbol{\scV}}}

\newcommand{\mt}{{\scM_\theta}}
\newcommand{\xfano}{{\bfX_\bSigmabar}}
\newcommand{\kfano}{{\bfX_\bSigma}}
\newcommand{\ksing}{{X_\Sigma}}
\newcommand{\kres}{{\widetilde{X_\Sigma}}}
\newcommand{\lres}{\scL'}

\newcommand{\lsing}{L}
\newcommand{\lk}{{\bscL}}
\newcommand{\vk}{{\bscV}}
\newcommand{\vfano}{{\bscE}}

\newcommand{\crep}{{X_{\Sigmatilde}}}
\newcommand{\vcrep}{{\scV}}
\newcommand{\lcrep}{{\scL}}

\newcommand{\dimfan}{n}
\newcommand{\nooc}{r}
\newcommand{\nooctilde}{{\widetilde{\nooc}}}

\title{Dimer models and exceptional collections}
\author{Akira Ishii and Kazushi Ueda}
\date{}
\begin{document}

\maketitle
\begin{abstract}
We construct a full strong exceptional collection
consisting of line bundles
on any two-dimensional smooth toric weak Fano stack.
The total endomorphism algebra of the resulting collection is isomorphic
to the path algebra of a quiver with relations
associated with a dimer model and a perfect matching on it.
%We also discuss the torus-equivariant version of the story.
\end{abstract}

\section{Introduction} \label{sc:introduction}
A dimer model is a bicolored graph on a real 2-torus
which encodes the information of a quiver with relations.
The main result of \cite{Ishii-Ueda_DMSMC} states that
for any smooth quasi-projective toric Calabi-Yau 3-fold $M$,
there is a dimer model $G$ such that
\begin{itemize}
 \item
the moduli space $\scM_\theta$ of $\theta$-stable representations
of the quiver $\Gamma$ with relations associated with $G$
of dimension vector $(1, \dots, 1)$
%which is stable with respect to the stability parameter $\theta$
%in the sense of King \cite{King}
is isomorphic to $M$ if we choose a suitable stability parameter $\theta$, and
 \item
the direct sum $\vcrep = \bigoplus_v \lcrep_v$
of the tautological bundles
is a tilting object
whose endomorphism algebra is isomorphic to the path algebra
$\bC \Gamma$ of the quiver $\Gamma$ with relations.
\end{itemize}
This gives a description of the derived category of coherent sheaves
on any smooth toric Calabi-Yau 3-fold
in terms of a quiver with relations;
$$
 D^b \coh M \cong D^b \module \bC \Gamma.
$$
The same result can also be obtained
by combining the existence of an isoradial dimer model
by Gulotta \cite{Gulotta},
the Calabi-Yau property of an isoradial dimer model
by Broomhead
\cite{Broomhead}
(cf. also \cite{Mozgovoy-Reineke,
Davison,
Bocklandt_CCDM,
Ishii-Ueda_CCDM})
and the Calabi-Yau trick by
Bridgeland, King and Reid \cite{Bridgeland-King-Reid}
(cf. also \cite{Van_den_Bergh_NCR}).

The aim of this paper is to give a similar description
for the derived category of coherent sheaves
on a two-dimensional smooth toric weak Fano stack.
Here, a smooth toric Deligne-Mumford stack
with the trivial generic stabilizer is said to be
a {\em weak Fano stack}
if the anti-canonical bundle is nef and big.
Let $\bfX$ be such a stack and
$K_\bfX$ be the total space of its canonical bundle.
%The stacky fan describing $\bfX$ is determined
%by the set $\{ v_i \}_{i=1}^\nooc$ of generators
%of the one-dimensional cones.
Then a projective crepant resolution
of the coarse moduli space of $K_\bfX$
is a smooth toric Calabi-Yau 3-fold $M$,
and one has an equivalence
$$
 \Phi : D^b \coh M \to D^b \coh K_\bfX
$$
of derived categories.
The image $\Phi(\vcrep)$ of the tilting object on $M$
is a tilting object on $K_\bfX$,
which gives a generator on $\bfX$
by the derived restriction to the zero-section.

%Despite the simplicity,
%this construction has the following drawbacks:
%\begin{itemize}
% \item
%The resulting generator on $\bfX$ is
%not necessarily a vector bundle, or even a sheaf,
%but only a complex of sheaves.
% \item
%It is not obvious if the generator is a tilting object.
%Moreover,
%the relation between the endomorphism algebra of the generator
%and the path algebra of the quiver with relations is unclear.
%\end{itemize}
%%
%To remedy these problems,
%we use the following alternative construction;
%%
The basic strategy is to find a suitable equivalence $\Phi$
so that the resulting generator on $\bfX$ will be
not only an object of the derived category
but a direct sum of line bundles.
Then its restriction to the zero section will be a tilting object.
To be more precise,
we first construct line bundles before choosing
a derived equivalence as follows:
\begin{enumerate}
 \item
We start with the moduli space $\mt$ with an arbitrary generic stability parameter $\theta$,
which may not be lying over the coarse moduli space of $K_\bfX$.
This will enable us to use an arbitrary {central perfect matching}
to describe the quiver with relations that is derived equivalent to
$\bfX$ in Theorem \ref{th:total_morphism_algebra}.
 \item
Let $\lk_v$ be the line bundle on $K_{\bfX}$ obtained as the
proper transform of $\lcrep_v$ on $\mt$ and
put $\vk = \bigoplus_v \lk_v$.
Then one has an isomorphism
$$
 \End(\vk) \cong \End(\vcrep) \cong \bC \Gamma.
$$
 \item
The acyclicity
$$
 \Ext^k(\vcrep, \vcrep) = 0, \qquad k \ne 0
$$
of $\vcrep$ implies that of $\vk$.
 \item
Now one can use
%\begin{itemize}
% \item
%the acyclicity of $\vk$,
% \item
%the finiteness of the homological dimension
%of $\End(\vk) \cong \bC \Gamma$,
% \item
%the indecomposability of $D^b \coh K_\bfX$, and
% \item
%the triviality of the relative Serre functor on $D^b \coh K_\bfX$
%\end{itemize}
the Calabi-Yau trick
\cite{Bridgeland-King-Reid, Van_den_Bergh_NCR,
Bezrukavnikov-Kaledin_2004}
to show that $\vk$ is a tilting object.
 \item
The restriction $\vfano$ of $\vk$ to the zero-section
is a tilting object on $\bfX$,
which is a direct sum of line bundles.
 \item
The endomorphism algebra of $\vfano$
is the quotient of the endomorphism algebra of $\vk$
by the ideal consisting of elements
vanishing on the zero section.
The isomorphism
$
 \End(\vk) \cong \End(\vcrep)
$
gives a description of this ideal
in terms of a perfect matching on the dimer model.
\end{enumerate}

This gives a proof of a particular case
of a conjecture of King,
together with a description of the total morphism algebra:

\begin{theorem} \label{th:main}
Any two-dimensional smooth toric weak Fano stack
has a full strong exceptional collection
consisting of line bundles,
such that the total morphism algebra of the collection
is isomorphic to the path algebra
of a quiver with relations
associated with a consistent dimer model
and a perfect matching on it.
\end{theorem}

The original conjecture of King
\cite[Conjecture 9.3]{King_tilt}
states that a smooth complete toric variety
has a full strong exceptional collection
consisting of line bundles.
This is shown to be false
by Hille and Perling \cite{Hille-Perling},
who subsequently gave a necessary and sufficient condition
for a smooth complete toric surface
to have such a collection \cite{Hille-Perling_ESISRS}.
Kawamata \cite{Kawamata_DCTV} shows that
a smooth projective toric stack has
a full exceptional collection consisting of sheaves.
Borisov and Hua \cite{Borisov-Hua} suggested
to extend the conjecture to stacks,
with an additional assumption that the toric stack be weak Fano,
and proved it for toric Fano stacks
of Picard number or dimension at most two.
This modified conjecture turned out to be false
by Efimov \cite{Efimov_MLECLB}.
A fine moduli interpretation of any smooth projective toric varieties,
which was one of the original motivations of King,
is obtained by Craw and Smith \cite{Craw-Smith}.
%See also Bergman and Proudfoot \cite{Bergman-Proudfoot}.
The concept of dimer models and the idea to use them
to construct full strong exceptional collections
on toric surfaces came from string theorists;
see e.g. \cite{Franco-Hanany-Martelli-Sparks-Vegh-Wecht_GTTGBT,
Franco-Hanany-Vegh-Wecht-Kennaway_BDQGT,
Franco-Vegh_MSGTDM,
Hanany-Herzog-Vegh_BTEC,
Hanany-Kennaway_DMTD} and references therein.

The organization of this paper is as follows:
We collect basic facts on line bundles on toric stacks
in Section \ref{sc:toric}, and
basic definitions on dimer models
in Section \ref{sc:dimer}.
The relation between exceptional collections,
tilting objects and derived equivalences
is summarized in Section \ref{sc:tilting}.
We recall the main result of \cite{Ishii-Ueda_DMSMC}
in Section \ref{sc:tilting_CY3}.
A tilting object on the total space of the canonical bundle
of a toric weak Fano stack will be constructed
in Section \ref{sc:tilting_kfano},
which will be restricted to the image of the zero-section
to produce a full strong exceptional collection
in Section \ref{sc:tilting_Fano}.
In Section \ref{sc:tilting_divisors},
we use the same idea as in Section \ref{sc:tilting_Fano}
to give a description of the derived category of coherent sheaves
on the union of toric divisors in $\scM_\theta$
in terms of a dimer model.
%In Section \ref{sc:cyclic},
%we discuss the inverse construction of a dimer model
%from an exceptional collection of line bundle,
%and show as an immediate consequence of a theorem of
%Bocklandt \cite{Bocklandt_CYAWQP} that
%a full exceptional collection of line bundle comes from a dimer model
%if and only if it is a {\em strong cyclic collection}.

{\bf Acknowledgment}:
We thank Osamu Iyama for valuable discussions.
In particular,
Remark \ref{rem:APR} is due to him.
We also thank Hokuto Uehara for suggesting improvements,
and Markus Perling
for pointing out the reference \cite{Hille-Perling_ESISRS}
and posing a question which lead to Remark \ref{rm:cyclic}.
A.~I. is supported by Grant-in-Aid for Scientific Research (No.18540034).
K.~U. is supported by Grant-in-Aid for Young Scientists (No.20740037).
A part of this work is done
while K.~U. is visiting the University of Oxford,
and he thanks the Mathematical Institute for hospitality
and Engineering and Physical Sciences Research Council
(EP/F055366/1) for financial support.

\section{Line bundles on toric stacks}
 \label{sc:toric}

We recall the definition of toric stacks from \cite{Borisov-Chen-Smith}.
Let $N$ be a free abelian group of rank $\dimfan$.
Note in \cite{Borisov-Chen-Smith} $N$ is allowed to have torsions
so that the associated stack may have generic stabilizers
but we consider the torsion free case in this paper.
A {\em stacky fan}
$
 \bSigma = (\Sigma, \{ v_i \}_{i=1}^\nooc)
$
consists of a fan $\Sigma$ in $N_\bR = N \otimes \bR$
and the set $\{ v_i \}_{i=1}^\nooc$ of generators
of one-dimensional cones in $\Sigma$.
The toric stack $\bfX_\bSigma$ associated with $\bSigma$ is defined
as the quotient stack
$$
 \bfX_\bSigma = [(\bC^\nooc \setminus \SR(\Sigma)) / K],
$$
where the Stanley-Reisner locus $\SR(\Sigma)$
consists of points $(z_1, \dots, z_\nooc)$ such that
there is no cone in $\Sigma$ which contains
all $v_i$ for which $z_i = 0$,
and
$$
 K = \Ker (\phi \otimes \bCx)
$$
is the kernel of the tensor product with $\bCx$
of the map
$
 \phi :  \Ntilde = \bZ^\nooc \to  N
$
sending the $i$-th coordinate vector $e_i$
to $v_i$
for $i = 1, \dots, \nooc$.
We sometimes write
$$
 U_\Sigma = \bC^\nooc \setminus \SR(\Sigma)
$$
so that
$$
 \bfX_\bSigma = [U_\Sigma / K].
$$
It follows from the definition that
the category of coherent sheaves on $X_\bSigma$ is equivalent
to the category of $K$-equivariant coherent sheaves on $U_\Sigma$;
$$
 \coh \bfX_\bSigma \cong \coh^K U_\Sigma.
$$
Let
$
 M = \Hom(N, \bZ)
$
and
$
 \Mtilde = \Hom(\Ntilde, \bZ)
$
be the abelian group dual to $N$ and $\Ntilde$ respectively.
The tori $\bT = \Spec \bC[M]$ and $\bTtilde = \Spec \bC[\Mtilde]$
act naturally on $X_\bSigma$ and $U_\Sigma$, and
the category of $\bT$-equivariant coherent sheaves
on $X_\bSigma$ is equivalent
to the category of $\bTtilde$-equivariant coherent sheaves
on $U_\Sigma$;
$$
 \coh^\bT \bfX_\bSigma \cong \coh^{\bTtilde} U_\Sigma.
$$
As a smooth toric variety,
$\Pic U_\Sigma$ is generated by invariant divisors
$D_i = \{ z_i = 0 \}$,
which are clearly trivial.
Hence one has
$$
 \Pic^\bT \bfX_\bSigma
  \cong \Hom(\bTtilde, \bCx)
  = \Mtilde.
$$
The Picard group of $X_\bSigma$ can be calculated
using the exact sequence
$$
 1 \to \Hom(\bT, \bCx) \to \Pic^\bT \bfX_\bSigma \to \Pic \bfX_\bSigma \to 1.
$$
The divisor $D_i = \{ z_i = 0 \}$
naturally corresponds
%to the character of $\bTtilde$ corresponding
to the $i$-th coordinate vector
in $\Mtilde$,
which defines a line bundle
$\scO(D_i) \in \Pic^\bT \bfX_\bSigma$.
Its image in $\Pic \bfX_\bSigma$ will again be denoted by $\scO(D_i)$.
Any line bundle $\scL \in \Pic^\bT \bfX_\bSigma$ can be represented as 
$
 \scO(D),
$
where $D = (f)$ is the divisor of
any $\bTtilde$-invariant rational section
$$
 f\in
   \lb
    \bC[z_1^{\pm 1}, \dots, z_\nooc^{\pm 1}] \otimes \scL
   \rb^{\bTtilde},
$$
which is unique up to scalar multiples.

The cohomology of $\bT$-equivariant line bundle is given as follows:

\begin{proposition} \label{prop:coh_of_toric_line_bundle}
Let $\bSigma$ be a simplicial stacky fan and
$\scO(D)$ be the $\bT$-equivariant line bundle
associated with a divisor
$
 D = \sum_{i=1}^\nooc a_i D_i.
$
Then the $\bT$-invariant part
of the cohomology group of $\scO(D)$ is given by
$$
 H^p_\bT(\bfX_\bSigma, \scO(D))
  \cong H_Z^p(|\Sigma|),
$$
where $|\Sigma|$ is the support of the fan $\Sigma$
underlying the stacky fan $\bSigma$,
$$
 \psi_D : |\Sigma| \to \bR
$$
is the piecewise-linear function
which is linear on each cone of $\Sigma$
satisfying
$$
 \psi_D(v_i) = a_i,
$$
and
$$
 Z = \{ x \in |\Sigma| \mid \psi_D(x) \ge 0 \}.
$$
\end{proposition}

See e.g. \cite[Section 3.5]{Fulton_ITV}
%or \cite[Proposition 4.1]{Borisov-Hua}
for a proof
of Proposition \ref{prop:coh_of_toric_line_bundle}.
The cohomology group of a line bundle
$\scL$ on $\bfX_\bSigma$ is the direct sum
$$
 H^p(\bfX_\bSigma, \scL)
  = \bigoplus_{D : \scL \cong \scO(D)} H^p_\bT(\bfX_\bSigma, \scO(D))
$$
over the set of toric divisors $D$ such that
$\scL \cong \scO(D)$,
which is a torsor over the lattice $M$.

\section{Dimer models} \label{sc:dimer}

%Dimer models are introduced by string theorists
%(see e.g. \cite{Franco-Hanany-Martelli-Sparks-Vegh-Wecht_GTTGBT,
%Franco-Hanany-Vegh-Wecht-Kennaway_BDQGT,
%Franco-Vegh_MSGTDM,
%Hanany-Herzog-Vegh_BTEC,
%Hanany-Kennaway_DMTD})
%to study supersymmetric quiver gauge theories
%in four dimensions.

A {\em dimer model} is a bicolored graph on a torus $T = \bR^2 / \bZ^2$
consisting of a set $B \subset T$ of black nodes,
another set $W \subset T$ of white nodes,
and a set $E$ of edges consisting of embedded line segments
connecting vertices of different colors.

A connected component of the complement
$T \setminus E$ is called a {\em face} of the graph.
A bicolored graph on $T$ is said to be a {\em dimer model}
if any face is simply-connected.
We only deal with dimer models
satisfying a consistency condition
described in \cite[Definition 3.5]{Ishii-Ueda_CCDM}.
See also \cite{Mozgovoy-Reineke, Davison, Broomhead, Bocklandt_CCDM}
for more about consistency conditions of dimer models.
Mathematical literature on dimer models also includes
\cite{Szendroi_NCDT,
Nagao-Nakajima,
Nagao_DCSCY,
Nagao_RONCDT,
Stienstra_Mahler_dimers,
Stienstra_dessins,
Stienstra_Chow,
Stienstra_CPA,
Mozgovoy,
Bender-Mozgovoy}.

A {\em quiver} consists of
a set $V$ of vertices,
a set $A$ of arrows, and
two maps $s, t: A \to V$ from $A$ to $V$.
For an arrow $a \in A$,
$s(a)$ and $t(a)$
are said to be the {\em source}
and the {\em target} of $a$
respectively.
A {\em path} on a quiver
is an ordered set of arrows
$(a_n, a_{n-1}, \dots, a_{1})$
such that $s(a_{i+1}) = t(a_i)$
for $i=1, \dots, n-1$.
We also allow for a path of length zero,
starting and ending at the same vertex.
The {\em path algebra} $\bC Q$
of a quiver $Q = (V, A, s, t)$
is the algebra
spanned by the set of paths
as a vector space,
and the multiplication is defined
by the concatenation of paths;
$$
 (b_m, \dots, b_1) \cdot (a_n, \dots, a_1)
  = \begin{cases}
     (b_m, \dots, b_1, a_n, \dots, a_1) & s(b_1) = t(a_n), \\
      0 & \text{otherwise}.
    \end{cases}
$$
A {\em quiver with relations}
is a pair of a quiver
and a two-sided ideal $\scI$
of its path algebra.
For a quiver $\Gamma = (Q, \scI)$
with relations,
its path algebra $\bC \Gamma$ is defined as
the quotient algebra $\bC Q / \scI$.

A dimer model $G = (B, W, E)$ encodes
the information of a quiver $\Gamma = (Q, \scI)$
with relations
in the following way:
The set $V$ of vertices of $Q$
is the set of faces of the graph, and
the set $A$ of arrows of $Q$
is the set $E$ of edges of the graph.
The directions of the arrows are determined
by the colors of the vertices of the graph,
so that the white vertex $w \in W$ is on the right
of the arrow.
In other words,
the quiver is the dual graph of the dimer model
equipped with an orientation given by
rotating the white-to-black flow on the edges of the dimer model
by minus 90 degrees.

For an arrow $a \in A$,
there exist two paths $p_+(a)$
and $p_-(a)$
from $t(a)$ to $s(a)$,
the former going around the white vertex
connected to $a \in E = A$ clockwise
and the latter going around the black vertex
connected to $a$ counterclockwise.
Then the ideal $\scI$
of the path algebra is
generated by $p_+(a) - p_-(a)$
for all $a \in A$.

A {\em perfect matching}
(or a {\em dimer configuration})
on a dimer model $G = (B, W, E)$
is a subset $D$ of $E$
such that for any vertex $v \in B \cup W$,
there is a unique edge $e \in D$
connected to $v$.
The two-sided ideal of $\bC \Gamma$
generated by arrows $a$ in $D \subset E = A$
will be denoted by $\scI_D$.
%
%The height function $h_{D, D_0}$ is
%a locally-constant function on
%$\bR^2 \setminus (D \cup D_0)$
%which increases (resp. decreases)
%by $1$
%when one crosses an edge $e \in D$
%with the black (resp. white) node
%on his right
%or an edge $e \in D_0$
%with the white (resp. black) node
%on his right.
%This rule determines the height function
%up to additions of constants.
%The height function may not be periodic
%even if $D$ and $D_0$ are periodic,
%and the {\em height change}
%$h(D, D_0) = (h_x(D, D_0), h_y(D, D_0)) \in \bZ^2$
%of $D$ with respect to $D_0$
%is defined as the difference
%\begin{align*}
% h_x(D, D_0) &= h_{D, D_0}(p+(1,0)) - h_{D, D_0}(p), \\
% h_y(D, D_0) &= h_{D, D_0}(p+(0,1)) - h_{D, D_0}(p) 
%\end{align*}
%of the height function,
%which does not depend on the choice of
%$p \in \bR^2 \setminus (D \cup D_0)$.
%More invariantly,
%height changes can be considered
%as an element of $H^1(T, \bZ)$.
%The dependence of the height change
%on the choice of the reference matching
%is given by
%$$
% h(D, D_1) = h(D, D_0) - h(D_1, D_0)
%$$
%for any three perfect matchings $D$, $D_0$ and $D_1$.
%We will often suppress the dependence of the height difference
%on the reference matching
%and just write $h(D) = h(D, D_0)$.

\section{Exceptional collections and tilting objects}
 \label{sc:tilting}

Let $\bfX$ be a smooth stack and
$
 \scT = D^b \coh \bfX
$
be the derived category of coherent sheaves on $\bfX$.

\begin{definition}
\ \vspace{-6mm} \\
\begin{enumerate}
 \item
An object $E$ of $\scT$ is {\em acyclic}
if $\Ext^k(E, E) = 0$ for $k \ne 0$.
 \item
An acyclic object $E$ is {\em exceptional}
if $\End(E)$ is spanned by the identity morphism.
 \item
A sequence $(E_1, \dots, E_n)$ of exceptional objects is
an {\em exceptional collection}
if
$
 \Ext^k(E_i, E_j) = 0
$
for $1 \le j < i \le n$.
 \item
An exceptional collection $(E_1, \dots, E_n)$ is {\em strong}
if $\Ext^k(E_i, E_j) \ne 0$ implies $k = 0$.
 \item
An exceptional collection is {\em full}
if it generates $\scT$
as a triangulated category.
 \item
An object $V$ is a {\em generator}
if $\RHom(V, X) = 0$ implies $X \cong 0$.
 \item
An acyclic generator
is called a {\em tilting object}.
\end{enumerate}
\end{definition}

%The following lemma is clear from the definition:
%\begin{lemma}
%\ \vspace{-6mm} \\
%\begin{enumerate}
% \item
%A line bundle on a smooth stack $\bfX$ is exceptional
%if and only if $\bfX$ is proper.
% \item
%A sequence $(E_1, \dots, E_n)$ of line bundles
%on a smooth proper stack
%is a full strong exceptional collection
%if and only if $\bigoplus_{i=1}^n E_i$ is a tilting object.
%\end{enumerate}
%\end{lemma}
%
Note that a sequence $(E_1, \dots, E_n)$ of line bundles
on a smooth proper stack
is a full strong exceptional collection
if and only if $V = \bigoplus_{i=1}^n E_i$ is a tilting object.
The algebra
$
 \End(V) = \bigoplus_{i, j=1}^n \Hom(E_i, E_j)
$
%of the tilting object
%obtained as the direct sum
%$V = \bigoplus_{i=1}^n E_i$
%of a full strong exceptional collection
%$(E_1, \dots, E_n)$
is called the {\em total endomorphism algebra} of the collection.
It is a finite-dimensional algebra
which can be described
as the path algebra of a quiver with relations.

A tilting object induces a derived equivalence;

\begin{theorem}[{\cite{Rickard, Bondal_RAACS}}]
If a smooth stack $\bfX$ has a tilting object $V$,
then the functor
$$
 \RHom(V, \bullet) :
  D^b \coh \bfX \to D^b \module \End(V)
$$
induces an equivalence of triangulated categories.
\end{theorem}

\section{A tilting bundle on a smooth toric Calabi-Yau 3-fold}
 \label{sc:tilting_CY3}

Let $\{ v_i \}_{i=1}^\nooc$ be a set of points
on a lattice $\Nbar$ of rank two, and
$
 \bSigmabar = ( \Sigmabar, \{ v_i \}_{i=1}^\nooc)
$
be a two-dimensional complete stacky fan
whose two-dimensional cones are
$$
 \sigmabar_i = \bR_+ v_i + \bR_+ v_{i+1}, \qquad i = 1, \dots, \nooc,
$$
with $v_{\nooc+1} = v_1$.
The toric stack $\xfano$ associated with $\bSigmabar$ is
a weak Fano stack if and only if all the $v_i$ are
on the boundary of the lattice polygon
$$
 \Delta = \Conv \{ v_i \}_{i=1}^\nooc
$$
defined as the convex hull of $\{ v_i \}_{i=1}^\nooc$.
The torus $\Spec \bC[\Mbar]$ acting on $X_\bSigmabar$
will be denoted by $\bTbar$,
where $\Mbar = \Hom(\Nbar, \bZ)$.
Let
$$
 p : \kfano \to \xfano
$$
be the total space of the canonical bundle of $\xfano$.
The stacky fan $\bSigma$
corresponding to the total space of the canonical bundle
%$\kfano$
is given by
$(\Sigma, \{ \vtilde_i \}_{i=0}^\nooc)$,
where the generators of one-dimensional cones are given by
$$
 \vtilde_0 = (0, 0, 1), \qquad
 \vtilde_i = (v_i, 1), \quad i = 1, \dots, \nooc
$$
and three-dimensional cones of $\Sigma$ consists of
$$
 \sigma_i = \bR_+ \vtilde_i + \bR_+ \vtilde_{i+1} + \bR_+ \vtilde_0,
  \qquad i = 1, \dots, \nooc.
$$
We write the lattice containing $\Sigma$
and its dual
as $N$ and $M$ respectively.
Let $\bT = \Spec \bC[M]$ be the torus
acting on $\kfano$.
The toric divisor associated with the one-dimensional cone
generated by $\vtilde_i$ will be denoted by $\bfD_i$.

Let $R = H^0(\scO_\kfano)$ be the coordinate ring of the three-dimensional affine toric variety
associated with the cone over $\Delta$.
The following are shown in \cite{Ishii-Ueda_08, Ishii-Ueda_DMSMC}:

\begin{itemize}
 \item
There is a consistent dimer model $G = (B, W, E)$
on $T = (\Mbar \otimes \bR) / \Mbar$
such that the moduli space $\mt$
of $\theta$-semi-stable representations of $\Gamma$
with dimension vector $(1, \dots, 1)$
is a crepant resolution of $\Spec R$
where $\theta$ is a generic stability parameter
for the quiver $\Gamma$ with relations
associated with the dimer model $G$.
 \item
For a prime toric divisor $D$ in $\mt$,
there is a perfect matching,
which we write $D$ again by abuse of notation,
such that the divisor is the zero locus
of the arrows contained in the perfect matching.
 \item
For any perfect matching $D$ on $G$,
there is a stability parameter $\theta$ such that
$D$ corresponds to a prime toric divisor on $\mt$.
 \item
The tautological bundle
$$
 \vcrep = \bigoplus_{v \in V} \lcrep_v
$$
on $\mt$ is a tilting object in $D^b \coh \mt$
such that
$$
 \End \vcrep \cong \bC \Gamma.
$$
\end{itemize}

The fan describing the toric variety $\mt$
is a refinement of the fan consisting of the cone over $\Delta$
and its faces, and will be denoted by $\Sigmatilde$.
The perfect matchings corresponding to $\vtilde_i$
for $i = 0, \dots, \nooc$ will be denoted by $D_i$.
The generators of the one-dimensional cones of $\Sigmatilde$
which does not belong to $\Sigma$
will be denoted by
$v_{\nooc + 1}, \ldots, v_{\nooctilde}$,
and the corresponding toric divisors will be written as
$D_{\nooc + 1}, \ldots, D_{\nooctilde}$.
A perfect matching which is not on a vertex of $\Delta$
depends on the choice of $\theta$
\cite[Proposition 6.5]{Ishii-Ueda_DMSMC}.
A perfect matching which corresponds to $\vtilde_0$
under some stability parameter $\theta$
is called a {\em central perfect matching}.

\section{A tilting bundle on the canonical bundle}
 \label{sc:tilting_kfano}

Let
$$
 \bvarphi : \kfano \to \ksing
$$
be the natural morphism from a stack to its coarse moduli space, and
$$
 \varphi : \kres \to \ksing
$$
be a crepant resolution.
Let further $\crep$ be an arbitrary crepant resolution of $\Spec R$,
so that $\kres$ and $\crep$ are isomorphic in codimension one.

For a line bundle $\lcrep$ on $\crep$,
let $\lres$ be its proper transform on
$
 \kres,
$
$$
 \lsing = (\varphi_* \lres)^{\vee \vee}
$$
be the reflexive sheaf of rank one on $X_\Sigma$
obtained as the double dual of the direct image of $\lres$, and
$$
 \lk = (\bvarphi^* \lsing)^{\vee \vee}
$$
be the line bundle on $\kfano$ obtained as the double dual
of the pull-back of $L$ to $\kfano$.
If $\lcrep$ is isomorphic to $\scO_{\crep}(\Dtilde)$
for a divisor $\Dtilde = \sum_{i = 0}^\nooctilde a_i D_i$,
then $\lsing$ is isomorphic to $\scO_{\ksing}(D)$
where $D = \sum_{i=0}^\nooc a_i D_i$ is obtained from $\Dtilde$
by forgetting the toric divisors contracted by $\varphi$,
and $\lk$ is isomorphic to $\scO_{\kfano}(\bfD)$
where $\bfD = \sum_{i=0}^\nooc a_i \bfD_i$ is the pull-back of
$D$ by $\bvarphi$.

\begin{lemma} \label{lem:isom_if_reflexive}
The spaces of global sections of $\lcrep$, $\lres$, $\lsing$ and $\lk$ are
related as follows:
$$
H^0(\crep, \lcrep) =
 H^0(\kres, \lres)
  \subset H^0(\ksing, \lsing)
  = H^0(\kfano, \lk).
$$
Moreover,
if $H^0(\crep, \lcrep)$ is a relfexive $R$-module,
then the inclusion in the middle is an isomorphism.
\end{lemma}
\begin{proof}
The first equality follows from the fact that
$\crep$ and $\kres$ are isomorphic in codimension one.
The inclusion
$
 \varphi_* \lres \subset \lsing
$
implies
$
 H^0(\kres, \lres)
  \subset H^0(\ksing, \lsing).
$
Since both $\lsing$ and $\bvarphi_* \lk$ are reflexive,
the inclusion
$
 \lsing \hookrightarrow \bvarphi_* \lk
$
is an isomorphism and
the second equality follows.
Finally, the four spaces have structures of torsion free $R$-modules which coincieds on the smooth locus of $\Spec R$.
Therefore, the reflexivity assumption implies the last assertion.
\end{proof}

The description of the cohomology
of a $\bT$-equivariant line bundle on $\crep$
admits the following simplification in the present situation:
Let
$
 D = \sum_{i=0}^\nooctilde a_i D_i
$
be a toric divisor on $\crep$
and
$
 \psi_D : |\Sigmatilde| \to \bR
$
be the piecewise-linear function
associated with $D$.
Put
$$
 \Delta = |\Sigmatilde| \cap \left(\bR^2 \times \{1\}\right)
$$
and
$$
 \Zbar = Z \cap \Delta,
$$
where
$$
 Z = \{ x \in |\Sigmatilde| \mid \psi_D(x) \ge 0 \}.
$$
Then one has
$$
 H^p_{\bT}(\kres, \scO(D))
  \cong H^p_{Z}(|\Sigma|)
  \cong H^p_{\Zbar}(\Delta)
$$
since everything is linear on the third coordinate.

The cohomology of a $\bT$-equivariant line bundle on $\kfano$
has an analogous description,
which can be simplified further as follows:
For a toric divisor
$
 \bfD = \sum_{i=0}^\nooc a_i \bfD_i
$
on $\kfano$,
let
$
 \sgn(\bfD)
  = (\sgn a_0; \sgn \bfDbar)
  = (\sgn a_0; \sgn a_1, \dots, \sgn a_\nooc)
$
be the sequence of the signatures of the coefficients of $\bfD$,
where
$$
 \sgn a =
  \begin{cases}
   + & a \ge 0, \\
   - & a < 0.
  \end{cases}
$$
Consider $\sgn(\bfDbar)$ as a cyclic sequence.
A {\em $-$-interval} in $\sgn(\bfDbar)$
is a succession of $-$ bounded by one $+$ and
the next $+$.
Then it follows
from Proposition \ref{prop:coh_of_toric_line_bundle}
that
\begin{align*}
 \rank H^0_\bT(\scO(\bfD)) &=
  \begin{cases}
   1 & \sgn(\bfD) = (+;+\cdots+), \\
   0 & \text{otherwise},
  \end{cases} \\
 \rank H^1_\bT(\scO(\bfD)) &=
  \begin{cases}
   0 & \sgn a_0 = - \text{ or } \sgn(\bfDbar) = (+;+\cdots+), \\
   \# \{-\text{-intervals}\} - 1 & \text{otherwise},
  \end{cases} \\
 \rank H^2_\bT(\scO(\bfD)) &=
  \begin{cases}
   1 & \sgn(\bfD) = (+;-\cdots-), \\
   0 & \text{otherwise}.
  \end{cases}
\end{align*}

\begin{proposition} \label{prop:acyclicity}
Assume that a line bundle $\lcrep$ and its dual $\lcrep^\vee$
on $\crep$ are acyclic.
\begin{enumerate}
\item $H^0(\lcrep)$ is Cohen-Macaulay and hence reflexive as an $R$-module.
\item
Assume further that $\lcrep$ satisfies
the following condition: \vspace{5mm}\\
\begin{tabular}{cc}
 $(\ast)$ & 
\begin{minipage}[]{.8 \linewidth}
if one has
$\sgn a_i = \sgn a_{j} = -$
where the line segment $[\vtilde_i, \vtilde_j]$
between $\vtilde_i$ and $\vtilde_j$
lies on the boundary of $\Delta$,
then $\sgn a_{k} = -$ for any $k$ such that
$\vtilde_k \in [\vtilde_i, \vtilde_j]$.
\end{minipage}
\end{tabular}
%\begin{list}{}{}
% \item[$\ast$]
%if one has
%$\sgn a_i = \sgn a_{j} = -$
%where the line segment $[\vtilde_i, \vtilde_j]$
%between $\vtilde_i$ and $\vtilde_j$
%lies on the boundary of $\Delta$,
%then $\sgn a_{k} = -$ for any $k$ such that
%$\vtilde_k \in [\vtilde_i, \vtilde_j]$.
%\end{list}
\vspace{5mm} \\
Then the corresponding line bundle $\lk$ on $\kfano$ is acyclic.
\end{enumerate}
\end{proposition}

\begin{proof}
The first statement follows from arguments
in \cite[Proposition A.2]{Toda-Uehara} as follows:
Since $R$ is Gorenstein and $f$ is crepant,
we have $f^{!} R \cong \scO_{\crep}$.
Then the Grothendieck duality for $f$
and the acyclicity of $\lcrep$ impliy
$$
 \Hom^i_R(H^0(\lcrep), R)
  \cong \Hom^i_{\crep}(\lcrep, f^{!}R)
  \cong \Hom^i_{\crep}(\lcrep, \scO_{\crep}),
$$
which is $0$ for $i \ne 0$ by the acyclicity of $\lcrep^{\vee}$.
This implies that $H^0(\lcrep)$ is a Cohen-Macaulay $R$-module,
and therefore is reflexive since $\dim R \ge 2$.

For the second statement,
first assume that $H^2_{\bT}(\scO(\bfD)) \ne 0$.
Then one has $\sgn(\bfD) = (+;-\cdots-)$
so that
$$
 \psi_{\Dtilde}|_\Delta(\vtilde) < 0
$$
for any $\vtilde \in \partial \Delta$.
It follows that
$
 \partial \Delta \subset \Delta \setminus \Zbar
$
and
$
  \vtilde_0 \in \Zbar,
$
which implies
$$
 H^2_\bT(\crep, \lcrep)
  \cong H^2_{\Zbar}(\Delta)
  \cong H^1(\Delta \setminus \Zbar)
%  \cong H^2_{\Zbar}(\Delta)
  \ne 0.
$$

Next assume that $H^1_\bT(\scO(\bfD)) \ne 0$,
which implies that
$\sgn a_0 = +$ and
there are more than two $-$-intervals
in $\sgn(\bfDbar)$.
If $\scO(\Dtilde)$ is acyclic,
then $\Delta \setminus \Zbar$ is connected (and simply-connected).
Now consider the toric divisor
$$
 \Dtilde^{\vee}
  = \sum_{i=0}^\nooctilde (- \atilde_i - 1) \Dtilde_i
$$
which satisfies
$$
 \scO(\Dtilde^{\vee}) \cong \scO(- \Dtilde)
$$
as a non-equivariant line bundle
(the subtraction of one from all the coefficients of $\Dtilde^\vee$
corresponds to a change of a $\bT$-linearization).
Put
\begin{equation*}
 \Zbar^\vee
  := \{ x \in \Delta \mid \psi_{\Dtilde^\vee}(x) \ge 0 \}%\\
   = \{ x \in \Delta \mid \psi_{\Dtilde}(x) \le -1 \}.
\end{equation*}
%Then $\Mbar \cap \Zbar^\vee$ is the complement of
%$\Mbar \cap \Zbar$ in $\Mbar \cap \Delta$.
Then $\Zbar$ is a deformation retract of $\Delta \setminus \Zbar^\vee$.
Now the condition that $\scL^\vee$ is acyclic
implies that $\Delta \setminus \Zbar^\vee$ and therefore $\Zbar$ are connected, which contradicts
the connectedness of
$\Delta \setminus \Zbar$ and
the existence of multiple connected components
of $\partial \Delta \cap (\Delta \setminus \Zbar)$.
\end{proof}

\begin{lemma} \label{lem:cond_ast}
Let $\lcrep_v$ and $\lcrep_w$ be the tautological bundles
on the moduli space $\mt$
of quiver representations
associated with a consistent dimer model
with generic stability parameter $\theta$ and put
$
 \lcrep = \lcrep_v^\vee \otimes \lcrep_w.
$
Then $\lcrep$ satisfies the condition $(\ast)$
in Proposition \ref{prop:acyclicity} where $\crep=\mt$.
\end{lemma}

\begin{proof}
It follows from \cite{Ishii-Ueda_DMSMC} that the union
$
 \bigcup_{D \in [\vtilde_i, \vtilde_j]} D
$
of perfect matchings on the line segment
$
 [\vtilde_i, \vtilde_j] \subset \partial \Delta
$
is the union of isolated edges with the zig-zag paths corresponding to these line segments
dividing the torus $\Mbar_\bR / \Mbar$ into strips,
and the open subscheme $U$ of the moduli space
$\scM_\theta \cong \crep$ consisting of quiver representations
such that any arrow not in
$
 \bigcup_{D \in [\vtilde_i, \vtilde_j]} D
$
is non-zero is isomorphic to the product
of the moduli space $\scN$ of representations of
the McKay quiver of type $A_n$ in dimension two
and a one-dimensional torus;
$$
 U \cong \scN \times \bCx.
$$
A tautological bundle $\lcrep_v$ on $\scM_\theta$ restricts
to the outer tensor product $\lcrep'_v \boxtimes \scO_\bCx$
of a tautological bundle $\lcrep'_v$
on $\scN$ and a trivial bundle on $\bCx$.
Since the line bundle $(\lcrep'_v)^\vee \otimes \lcrep'_w$
on $\scN$ is acyclic,
the proof of Lemma \ref{lem:cond_ast} reduces to
Lemma \ref{lem:acyclicity_McKay} below.
\end{proof}

\begin{lemma} \label{lem:acyclicity_McKay}
Let $\Sigma$ be a two-dimensional fan
whose two-dimensional cones are given by
$$
 \bR_+ v_i + \bR_+ v_{i+1}, \qquad i = 0, 1, \dots, n
$$
where
$$
 v_i = (i, 1) \in N \cong \bZ^2.
$$
Let
$$
 D = \sum_{i=0}^n a_i D_i
$$
be a divisor on $X_\Sigma$,
where $D_i$ is the toric divisor
associated with the one-dimensional cone of $\Sigma$
generated by $v_i$.
If a line bundle $\scO(D)$ is acyclic and
the signatures of $a_0$ and $a_n$ are negative,
$\sgn(a_0) = \sgn(a_n) = -$,
then the signatures of all the $a_i$ are negative.
\end{lemma}

\begin{proof}
This is a corollary of the following fact,
which in turn follows immediately
from Proposition \ref{prop:coh_of_toric_line_bundle}:
$H^2_{\bT}(\scO(D)) = 0$ for any divisor $D$,
and $H^1_{\bT}(\scO(D))$ is the number of $-$-intervals
minus one if there is any,
and zero otherwise.
\end{proof}

Now by applying Lemma \ref{lem:isom_if_reflexive} and
Proposition \ref{prop:acyclicity}
to $\lcrep_v^\vee \otimes \lcrep_w$
for tautological bundles $\lcrep_v$ and $\lcrep_w$,
one shows that
$\vk = \bigoplus_v \lk_v$ is an acyclic bundle
satisfying
$$
 \End(\vk) \cong \bC \Gamma.
$$

The following definition is due to Bezrukavnikov and Kaledin:

\begin{definition}[{%Bezrukavnikov and Kaledin
\cite[Definition 2.1]{Bezrukavnikov-Kaledin_2004}}]
A non-zero object of an abelian category is {\em almost exceptional}
if $\Ext^i(M, M) = 0$ for $i > 0$ and
the algebra $\End(M)$ has finite homological dimension.
\end{definition}

The equivalence
$$
 D^b \coh \scM_\theta \cong D^b \module \bC \Gamma.
$$
implies that $\bC \Gamma$ has finite homological dimension.
Then the acyclicity of $\vk$ and the isomorphism
$$
 \End(\vk) \cong \bC \Gamma
$$
shows the following:
\begin{lemma}
The vector bundle $\vk$ on $\kfano$ is almost exceptional.
\end{lemma}

The proof of \cite[Lemma 4.2]{Bridgeland-King-Reid} actually shows
the following slightly stronger statement:

\begin{lemma}
The derived category of coherent sheaves
on a smooth Deligne-Mumford stack without a generic stabilizer
is indecomposable.
\end{lemma}

There is a morphism
$$
 \pi : \kfano \to \Spec R.
$$
%Then $D^b \coh \kfano$ is an $R$-linear triangulated category.
Since $R$ is Gorenstein, $\kfano$ is smooth and
the morphism $\pi$ is crepant,
the Grothendieck duality implies that
the identity functor is a Serre functor
of $D^b \coh \kfano$ with respect to $R$,
in the sense that there is a functorial isomorphism
$$
 \bR \Hom_R(\bR \pi_* \bR \scHom_{\kfano}(\scF, \scG), R)
  \cong \bR \pi_* \bR \scHom_{\kfano}(\scG, \scF)
$$
satisfying the compatibility conditions
in \cite{Bondal-Kapranov_Serre}.
This suffices to show the following:

\begin{theorem} \label{th:derived_equiv}
The functor
$$
 \Phi = \RHom(\vk, \bullet) :
  D^b \coh \kfano \to D^b \module \bC \Gamma
$$
is an equivalence of triangulated categories.
\end{theorem}

\begin{proof}
The proof is completely parallel to
\cite[Proposition 2.2]{Bezrukavnikov-Kaledin_2004}:
The functor $\Phi$ has a left adjoint
$$
 \Psi = \bullet \Lotimes_{\bC \Gamma} \vk
  : D^b \module \bC \Gamma \to D^b \coh \kfano,
$$
which produces a semiorthogonal decomposition
\begin{equation} \label{eq:semiorth_decomp}
 D^b \coh \kfano
  = (C, C^{\bot}),
\end{equation}
where $C$ is the essential image of $\Psi$ and
$C^{\bot}$ is the right orthogonal of $\vk$.
The $R$-Calabi-Yau property of $D^b \coh \kfano$ implies
that \eqref{eq:semiorth_decomp} is an orthogonal decomposition,
and the indecomposability of $D^b \coh \kfano$ shows that
$C^\bot$ is empty.
\end{proof}

\section{A tilting bundle on a toric weak Fano surface}
 \label{sc:tilting_Fano}

We use the same notation as in the previous sections.

\begin{lemma} \label{lem:restriction}
If
$
 \scV
$
be a tilting object on $\kfano$
which is a direct sum of line bundles,
then the restriction
$$
 \scE = \iota^* \scV
$$
of $\scV$
by the zero-section
$$
 \iota : \xfano \to \kfano
$$
is again a tilting object.
\end{lemma}

\begin{proof}
$\scE$ is a generator since it is the restriction of a generator
to a closed subscheme.
%Note that $\scL = p^* \scE$.
Vanishing of higher $\Ext$-groups $\Hom^{> 0}(\scE, \scE)$
follows from that of $\Hom^{> 0}(\scV, \scV)$
by
\begin{align*}
 \Hom^* (\scV, \scV)
  &= \Hom^* (p^* \scE, p^* \scE) \\
  &= \Hom^* (\scE, p_* p^* \scE) \\
  &= \Hom^* (\scE, \scE \otimes p_* \scO_{\kfano}) \\
  &= \Hom^* (\scE, \scE \otimes \bigoplus_{n=0}^\infty
                    \scK_{\xfano}^{\otimes(- n)} ),
\end{align*}
where $\scK_{\xfano}$ is the canonical sheaf of $\xfano$ and
the isomorphism
$
 \scV \cong p^* \scE
$
comes from the assumption that $\scV$
is a direct sum of line bundles.
\end{proof}

Note that the existence of a tilting object in $D^b \coh \xfano$
which is a direct sum of line bundles
is equivalent to the existence of a full strong exceptional collection
consisting of line bundles.

\begin{theorem} \label{th:total_morphism_algebra}
Let $D_0$ be an arbitrary central perfect matching.
Then there is a full strong exceptional collection
on $\xfano$ consisting of line bundles
such that the total morphism algebra satisfies
$$
 \End(\vfano) \cong \bC \Gamma / \scI_{D_0}.
$$
\end{theorem}

\begin{proof}
Choose a generic stability parameter $\theta$
such that $D_0$ is $\theta$-stable
as in \cite[Lemma 6.2]{Ishii-Ueda_08}.
Then we can apply Theorem \ref{th:derived_equiv}
and Lemma \ref{lem:restriction}
to see that the restrictions of $\lk_v$
form an full strong exceptional collection.

To obtain the description of the endomorphism algebra
$\End(\vfano)$ of $\vfano$, first note that
$\End(\vfano)$ is
the quotient of $\End(\vk)$ by the ideal
generated by sections vanishing at the toric divisor $\bfD_0$,
which is the image of the zero-section
$
 \iota : \xfano \to \kfano.
$
Now the theorem follows from the following facts,
which are obvious:

\begin{enumerate}
 \item
An element of $\End(\vk)$ vanishes
on the divisor $\bfD_0$ in $\kfano$
if and only if it vanishes at the generic point of $\bfD_0$.
 \item
This is equivalent to the vanishing of the corresponding element
of $\End(\vcrep)$ at the generic point of $D_0 \subset \mt$.
 \item
A path of the quiver gives an element of $\End(\vcrep)$
vanishing on $D_0$
if and only if it is contained in $\scI_{D_0}$.
\end{enumerate}
\end{proof}

\begin{remark} \label{rem:APR}
It follows from Theorem \ref{th:total_morphism_algebra}
that the derived category
$D^b \module \bC \Gamma / \scI_{D_0}$
of modules over $\bC \Gamma / \scI_{D_0}$
does not depend on the choice of a central perfect matching $D_0$.
Osamu Iyama pointed out that
this result also follows from the theory of
{\em 2-APR tilting}
\cite{Iyama-Opppermann_nRFA}.
%\cite{Auslander-Platzeck-Reiten, Iyama-Opppermann_nRFA}.
\end{remark}

%One can strengthen Theorem \ref{th:total_morphism_algebra} as follows:
%
%\begin{theorem} \label{th:any_central_pm}
%For any perfect matching $D$
%whose height change with respect to the central perfect matching
%$D_0$ is zero,
%there is a full strong exceptional collection on $\xfano$
%consisting of line bundles such that
%the total endomorphism algebra is isomorphic to
%$\bC \Gamma / \scI_D$.
%\end{theorem}
%
%\begin{proof}
%Choose a stability parameter $\theta$
%such that $\scM_\theta$ is smooth and
%$D$ is the central perfect matching
%\cite[Lemma 6.2]{Ishii-Ueda_08}.
%Although $\scM_\theta$ may not be a resolution of $\ksing$,
%one can perform a sequence of flops
%to obtain a projective crepant resolution $\kres$ of $\ksing$.
%The strict transform of the tautological bundle
%$\vcrep = \bigoplus_v \lres_v$ on $\scM_\theta$
%produces an acyclic bundle $\vcrep'$ on $\kres$,
%which turns out to be a tilting object by the Calabi-Yau trick
%since the endomorphism ring is invariant
%under flops:
%$$
% \End(\vcrep') \cong \End(\vcrep) \cong \bC \Gamma.
%$$
%Now one can produce a tilting object $\vk$ on $\kfano$ such that
%$
% \End(\vk) \cong \End(\vcrep')
%$
%just as in Section \ref{sc:tilting_Fano},
%and the ideal of $\End(\vk)$
%consisting of endomorphisms vanishing on the image of the zero section
%$
% \iota : \xfano \hookrightarrow \kfano
%$
%corresponds to the ideal of $\bC \Gamma$
%generated by arrows in $D$.
%\end{proof}

\begin{remark} \label{rm:cyclic}
The collection of line bundles obtained in Theorem
\ref{th:total_morphism_algebra} is
not only a full strong exceptional collection
but satisfies the condition that
the rolled-up helix algebra
$$
 \scA = \bigoplus_{i, j=1}^n \bigoplus_{k, n=0}^\infty
  \Hom^k(E_i, E_j \otimes \scK_{\xfano}^{\otimes (-n)})
$$
is concentrated in $k = 0$.
The converse statement that
any collection of line bundles satisfying this condition
comes from a dimer model is an immediate consequence
of the main result of Bocklandt \cite{Bocklandt_CYAWQP}
(cf. also \cite[Theorem 3.7]{Bocklandt_GTNCR}).
Indeed,
given such a collection $(E_1, \ldots, E_n)$,
one has
$
 \scA \cong \End_{\kfano}(\oplus p^* E_i)
$
where $p : \kfano \to \xfano$ is the canonical bundle of $\xfano$.
Then it is a toric non-commutative crepant resolution of
the coordinate ring $R = H^0(\scO_{\kfano})$
of the affine toric 3-fold
obtained by contracting the zero-section
by \cite[Proposition A.2]{Toda-Uehara}.
\end{remark}

\section{A tilting bundle on the union of toric divisors}
 \label{sc:tilting_divisors}

Let $Y = \scM_\theta$ be the smooth toric Calabi-Yau 3-fold
obtained as the moduli space of representations
of the quiver with relations
associated with a consisteint dimer model $G = (B, W, E)$ and
$\scV$ be the tilting object on $Y$
obtained as the direct sum of tautological line bundles
as in Section \ref{sc:introduction}.
%The restriction of the tautological bundle
%$\scV$ on $Y$ to $Y_0$ will be denoted by $\scV_0$.
For a vertex $v \in V$ of the quiver $(V, A, s, t)$
associated with $G$,
the {\em small cycle} $\omega_v \in \bC \Gamma$
is defined as $\omega_v = p_+(a) \cdot a$,
where $a$ is any arrow such that $s(a) = v$.
Let further
$
 W = \sum_{v \in V} \omega_v
$
be the central element of $\scA \cong \bC \Gamma$
obtained as the sum of the small cycles $\omega_v$
starting from each vertex $v \in Q_0$,
and $\scA_0 = \scA / (W)$ be the quotient ring
by the two-sided ideal generated by $W$.
Since the center of $\scA \cong \End(\scV)$
is isomorphic to $R \cong H^0(\scO_Y)$,
the element $W$ defines a regular functionon on $Y$.
Let $\iota : Y_0 \hookrightarrow Y$ be the inclusion
of the zero locus of $W$,
which is the union of toric divisors.

\begin{lemma}
The restriction $\scV_0 = \iota^* \scV$
is a tilting object in $D^b \coh Y_0$,
whose endomorphism algebra is isomorphic to $\scA_0$.
\end{lemma}

\begin{proof}
$\scV_0$ is a generator
just as in Lemma \ref{lem:restriction}
%since it is the restriction of a generator to a closed subscheme.
One has
\begin{align*}
 \RHom_{Y_0}(\scV_0, \scV_0)
  &= \RGamma(\scV_0^\vee \otimes \scV_0) \\
  &= \RGamma(\iota^* (\scV^\vee \otimes \scV)) \\
  &= \RGamma(\iota_* \iota^* (\scV^\vee \otimes \scV)) \\
  &= \RGamma(\{ \scV^\vee \otimes \scV
             \xto{W} \scV^\vee \otimes \scV \}) \\
  &= \{ \End(\scV) \xto{W} \End(\scV) \}
\end{align*}
which shows that $\End(\scV_0) \cong \End(\scV) / W \End(\scV)$ and
$\Ext_{Y_0}^i(\scV_0, \scV_0) = 0$ for $i > 0$.
%If a coherent sheaf $\scE$ on $Y_0$ satisfies
%$\RHom(\iota^* \scV, \scE) = 0$,
%then one has $\RHom(\scV, \iota_* \scE) = 0$ by adjunction,
%which implies $\iota_* \scE \cong 0$
%since $\scV$ is a generator.
%Since $\iota$ is a closed embedding,
%$\iota_* \scE \cong 0$ implies $\scE \cong 0$
%so that $\scV_0$ is a generator.
\end{proof}

Since a tilting object induces an equivalence
of bounded derived categories
(see e.g. \cite[Lemma 3.3]{Toda-Uehara}
for a proof without smoothness assumption),
we have the following:

\begin{corollary} \label{cor:dbcoh}
The functor
$$
 \Phi_0 = \RHom(\scV_0, \bullet)
  : D^b \coh Y_0 \to D^b \module \scA_0
$$
is an equivalence of triangulated categories.
\end{corollary}

The {\em bounded stable derived category} of $Y_0$ is the quotient category
$$
 \Dbsing(Y_0) = D^b \coh Y_0 / \perf Y_0.
$$
of the bounded derived category of coherent sheaves on $Y_0$
by the full subcategory consisting of perfect complexes
(i.e. bounded complexes of locally-free sheaves).
The bounded stable derived category of $\scA_0$ is defined similarly
as the quotient category
$$
 \Dbsing(\scA_0) = D^b \module \scA_0 / \perf \scA_0
$$
of the bounded derived category
$D^b \module \scA_0$ of finitely-generated right $\scA_0$-modules
by the full triangulated subcategory $\perf \scA_0$
consisting of perfect complexes
(i.e. bounded complexes of projective modules).
Since perfect complexes are characterized
in a purely categorical way
as {\em homologically finite} objects
\cite{Orlov_TCSEBLGM}
(i.e. objects $A$ such that
for any object $B$,
the group $\Hom(A, B[i])$
is trivial for all but a finite number of $i \in \bZ$),
an equivalence of bounded derived categories
induces an equivalence of bounded stable derived categories:

\begin{corollary} \label{cor:dbsing}
One has an equivalence
$$
 \Dbsing(Y_0) \cong \Dbsing(\scA_0)
$$
of %$\bZ / 2 \bZ$-graded 
triangulated categories.
\end{corollary}

Stable derived categories are introduced by Buchweitz
\cite{Buchweitz_MCM}
motivated by the theory of {\em matrix factorizations}
by Eisenbud \cite{Eisenbud_HACI}.
They are rediscovered by Orlov \cite{Orlov_TCS}
under the name
{\em triangulated categories of singularities}
following an idea of Kontsevich,
and plays essential role in homological mirror symmetry.
In particular,
the stable derived category $\Dbsing(Y_0)$
is expected to be equivalent to the derived category
of the {\em wrapped Fukaya category} of an affine curve
(see \cite{Abouzaid-Auroux-Efimov-Katzarkov-Orlov}
and references therein).
On the other hand,
the stable derived category
can be regarded as the derived category of a curved algebra
(i.e. a pair of an algebra and its central element,
cf. e.g. \cite{Positselski_TKDC}),
and the curved algebra $(\bC \Gamma, W)$
is the colimit of the following covariant functor $\scF$
from the category $\scC$ to the category of curved algebras:
A dimer model gives a one-dimensional CW complex
with nodes as 0-cells and edges as 1-cells,
and the category $\scC$ is the category
whose objects are open stars of 0-cells and 1-cells and
whose morphisms are inclusions.
In other words,
the category $\scC$ has nodes and edges as objects,
and there is a unique non-identity morphism
$e \to n$ for each adjacency of an edge $e \in E$
and a node $n \in B \sqcup W$.
The functor $\scF$ sends an edge $e \in E$
to the curved algebra $\scF(e) = (\scA_e, W_e)$
consisting of the path algebra $\scA_e$
of the cyclic quiver
$
 s_e \xto{a_e} t_e \xto{p_e} s_e
$
with two vertices $s_e$, $t_e$ and
two arrows $a_e$, $p_e$
and the central element
$W_e = a_e \cdot p_e + p_e \cdot a_e$.
For a node $n$,
let $(e_0, \ldots, e_r)$ be the set of edges
adjacent to $n$,
ordered clockwise if $n \in B$ and
counter-clockwise if $n \in W$.
%(only the cyclic ordering matters).
Then the value $\scF(n)$ for a node $n$ is
the curved algebra $(\scA_n, W_n)$
consisting of the path algebra $\scA_n$
of the cyclic quiver
$$
 s_{e_0}
  \xto{a_{e_0}} t_{e_0} = s_{e_1}
  \xto{a_{e_1}} t_{e_1} = s_{e_2}
  \xto{a_{e_2}} \cdots
  \xto{a_{e_r}} t_{e_r} = s_{e_0}
$$
with $r+1$ vertices and $r+1$ arrows
and the central element
$$
 W_n = \sum_{i=0}^r a_{e_{i-1}} \cdots a_{e_0}
      a_{e_r} \cdots a_{e_{i+1}} a_{e_i}.
$$
For each adjacency $e_i \to n$,
the map $\scF(e_i \to n) : \scF(e_i) \to \scF(n)$ sends
$
 a_{e_i}
$
and
$
 p_{e_i}
$
in $\scA_{e_i}$
to
$
 a_{e_i}
$
and
$
 a_{e_{i-1}} \cdots a_{e_0} a_{e_r} \cdots a_{e_{i+1}}
$
in $\scA_n$ respectively.
It is an interesting problem to relate this
to an idea of Kontsevich
\cite{Kontsevich_SGHA}
to describe the Fukaya category of a Stein manifold
in terms of a constructible sheaf of dg categories
on its Lagrangian skeleton.

\bibliographystyle{amsalpha}
\bibliography{bibs.bib}

\noindent
Akira Ishii

Department of Mathematics,
Graduate School of Science,
Hiroshima University,
1-3-1 Kagamiyama,
Higashi-Hiroshima,
739-8526,
Japan

{\em e-mail address}\ : \ akira@math.sci.hiroshima-u.ac.jp

\ \\

\noindent
Kazushi Ueda

Department of Mathematics,
Graduate School of Science,
Osaka University,\\
Machikaneyama 1-1,
Toyonaka,
Osaka,
560-0043,
Japan.

{\em e-mail address}\ : \  kazushi@math.sci.osaka-u.ac.jp
\ \vspace{0mm} \\

\end{document}